\documentclass[letterpaper,10pt]{IEEEtran}
\usepackage{cite}
\usepackage{graphics}
\usepackage{graphicx}
\usepackage{amsfonts,amsmath,amssymb}
\usepackage{enumerate}
\usepackage{hyperref}

\title{\huge On the Benefit of Nonlinear Control for Robust
Logarithmic\\ Growth: Coin Flipping Games as a Demonstration Case  }

\author{\Large Anton V. Proskurnikov and B. Ross Barmish
\thanks{Anton V. Proskurnikov is associate professor with Department of Electronics and Telecommunications at Politecnico di Torino, Turin, Italy and B. Ross Barmish is an emeritus professor in ECE at University of Wisconsin, Madison, and CEO of Robust Trading Solutions, LLC.}
\thanks{Emails: \texttt{anton.p.1982@ieee.org, bob.barmish@gmail.com}}
}

\def\be{\begin{equation}}
\def\ee{\end{equation}}
\def\ben{\begin{equation*}}
\def\een{\end{equation*}}

\def\boldP{\mathbb{P}}
\def\boldE{\mathbb{E}}

\def\ELG{\mathop{\sf ELG}\nolimits}
\def\IELG{\mathop{\sf IELG}\nolimits}

\begin{document}

\vskip -1in
\maketitle

\begin{abstract}
The takeoff point for this paper is the voluminous body of literature addressing recursive betting games with  expected logarithmic growth of wealth being the performance criterion. Whereas almost all existing papers involve use of linear feedback, the use of nonlinear control is conspicuously absent. This is epitomized by the large subset of this literature dealing with Kelly Betting. With this as the high-level motivation, we study the potential for use of nonlinear control in this framework. To this end, we consider a ``demonstration case'' which is one of the simplest scenarios encountered in this line of research: repeated flips of a biased coin with probability of heads~$p$, and even-money payoff on each flip. First, we formulate a new robust nonlinear control problem which we believe is both simple to understand and apropos for dealing with concerns about distributional robustness; i.e., instead of assuming that~$p$ is perfectly known as in the case of the classical Kelly formulation, we begin with a bounding set ~${\cal P} \subseteq [0,1]$ for this probability. Then, we provide a theorem, our main result, which gives a closed-form description of the optimal robust nonlinear controller and a corollary which establishes that it robustly outperforms linear controllers such as those found in the literature. A second, less significant, contribution of this paper bears upon the computability of our solution. For an $n$-flip game, whereas an admissible controller has~$2^n-1$ parameters, at the optimum only~$O(n^2)$ of them turn out to be distinct.
Finally, it is noted that the initial assumptions on payoffs
and the use of the uniform distribution on~$p$ are made solely for simplicity of the exposition and compliance with length requirements for a Letter. Accordingly, the paper also includes a new section with a discussion indicating how these assumptions can be relaxed.
\end{abstract}
\begin{IEEEkeywords}
Robust Control, Finance, Markov Processes
\end{IEEEkeywords}

\section{Introduction}
This paper addresses a large class of betting games described by discrete-time Markov processes. In this setting, the bettor
begins with initial account value~$V_0 > 0$.  At each stage~$k$, the control~$u_k$, alternatively called
the {\it betting strategy}, determines the size of the~$k$-th wager. Then, over~$n$ steps, assuming independent and identically
distributed random variables~$X_k$ as the {\it returns} with a known probability distribution, the resulting account value
trajectory~$V_1,V_2,...,V_n$, emanating from initial condition~$V_0 > 0$, is obtained recursively by
$$
V_{k+1} = V_k + u_kX_k.
$$
Within this context, the literature most closely related to this paper concentrates
on the design of a causal controller~$u$ maximizing the resulting  {\it Expected Logarithmic Growth}~(ELG)
$$
\ELG_u = \frac{1}{n} \boldE\log\left(\frac{V_n}{V_0}\right)
$$
subject to {\it budget constraints}~$|u_k| \leq V_k$ for~$k = 0,1,...,n-1$.

Perhaps, the most celebrated work along the lines above is the seminal paper by Kelly~\cite{Kelly_1956}; see also the early recognition of the power of the ELG approach  in~\cite{Breiman_1961,Thorp_1961,Thorp_1969,Cover_1984}. Over the decades  to follow, we see a voluminous body of literature, comprised of hundreds of papers, dealing with applications, extensions, and generalizations of Kelly's result in various directions.  We also see many papers providing rationale for the use of the logarithmic growth criterion versus other performance metrics. A selection of highlights from this work includes the detailed coverage of these topics in textbooks such as~\cite{Cover_Thomas_1991} and~\cite{Luenberger_1998} and the extensive collection of papers in~\cite{Maclean_et_al_2011}. It is also important to point out that this body of literature being cited includes major results on various properties of the ELG-maximizing controller over and above optimal logarithmic growth. That is, many authors provide results bearing on the ``asymptotic superiority'' of the ELG maximizer and cover other topics such as the relaxation the~i.i.d. assumption on the~$X_k$. In this regard, some good starting points for the uninitiated reader are~\cite{Hakansson_1971,Algoet_and_Cover_1988,Obrien_et_al_2021}. Finally, we draw attention to the doctoral dissertation of Hsieh~\cite{Hsieh_2019} which includes not only a comprehensive review of the earlier literature but also details and citations of his contributions and those of others over the preceding years.

To complete this brief perspective of the literature related to this paper, it is also important to mention the body of work dealing with ``distributional robustness'' issues arising in stochastic optimization; see\cite{Barmish_and_Shcherbakov_1999} and~\cite{Lagoa_and_Barmish_2002} where this terminology is introduced, the development of the theory in~ \cite{Delage_and_Ye_2010}, the important 2016 paper~\cite{Rujeerapaiboon_et_al_2016}, dealing specifically with distributional robustness in an ELG context and the more recent ongoing work~\cite{Li_2023} for along these lines.

Given the research setting above, the primary motivation for this paper is the fact that in the existing ELG literature, a problem formulation as one of nonlinear control is conspicuously absent; i.e., only linear control is considered. Whereas it is arguable, based on some of the results in existing work, that a nonlinear control cannot outperform a linear feedback when the probability distribution for the~$X_k$ is perfectly known, our main contention in this paper is that the same does not hold true when uncertainty in the underlying probability distributions is in play. Said another way, our main results provide compelling evidence that there are a large number of scenarios, involving distributional robustness considerations for which a nonlinear controller can outperform the ``best'' linear controller; e.g., see~\cite{Rujeerapaiboon_et_al_2016} and~\cite{Sun_and_Boyd_2018}. To this end, our analysis to follow demonstrates the potential of robust nonlinear control by considering one of the simplest possible ELG scenarios: a coin-flipping game with uncertainty in the probability of heads~$p$. Instead of taking~$p$ to be perfectly known as in the case of the classical Kelly formulation, we begin with  a bounding set~${\cal P} \subseteq [0,1]$ for this  probability. In this setting, the main result in this paper is a theorem which provides a complete closed-form solution of an optimal robust nonlinear control problem.
%Furthermore, to support the claim in the title that there is a ``benefit'' associated with nonlinear control,
As a corollary, we prove that our nonlinear controller robustly outperforms any linear controller.

A second less significant contribution of this paper bears on the computability of our new solution.
Whereas an admissible controller has~$2^n-1$ design parameters associated with the sample path points for an $n$-flip game,
surprisingly, at the optimum, many of them turn out to be the same with the resulting number of ``free parameters'' being of~$O(n^2)$.
Finally, we provide some initial illustrations bearing on the ``cost of imprecision'' in our knowledge of the probability~$p$.
To this end, some comparisons are made between the expected logarithmic growth associated with the optimal robust nonlinear
controller and the so-called {\it perfect-information Kelly optimum}.

\subsection*{Our Demonstration Case: Coin-Flipping}

To demonstrate the potential for the use of nonlinear control, we analyze
one of the simplest and most fundamental problems in the logarithmic growth literature: Making bets on~$n$ consecutive flips of a biased coin with probability of heads being~$p$ and even-money payoff. This initial assumption on the payoffs is made solely for simplicity brevity of the exposition; see Section~VII for a generalization to the unequal payoff case. Our simplified framework enables us to explain the key ideas behind our new nonlinear control formulation without being encumbered by additional
technical details.

Indeed, to demonstrate the potential for consideration of nonlinear controllers in future research, we begin with the widely celebrated betting scheme of Kelly~\cite{Kelly_1956}. That is, using the notation above, at stage~$k$, the controller generates the bet size as a linear feedback. That is,~$u_k= KV_k$ with  the understanding that~$u_k > 0$ and $u_k < 0$  corresponds to bets on heads and tails respectively. Consistent with this, even money payoffs are~$X_k = 1$ for heads and $X_k = -1$ for tails. 

Then, with probability of heads~$p$ assumed to be perfectly known and budget constraint $|u(k)| \leq V_k$ imposed, a straightforward calculation leads to  Kelly's optimal ELG maximizing feedback gain~$K^* = 2p-1$. It is important to note that, this simple linear controller proves to be optimal in many cases other than the simple scenario described above; e.g., in the widely cited 1971 paper by Hakansson~\cite{Hakansson_1971} (see also~\cite{Cover_Thomas_1991}), the bet size~$u_k$ may depend on the entirety of the past history~$V_0,V_1,...,V_{k-1},V_k$. However, as indicated earlier, our goal is to demonstrate the importance on nonlinear control when~$p$ is imperfectly known with robustness being a concern.

%\hskip -.5in
\section{Admissible Nonlinear Controllers and\\
 Resulting Expected Logarithmic Growth}

For the coin-flipping game at hand with sample path space
$$
{\cal X} \doteq \{-1,1\}^n,
$$
a mapping~$K: {\cal X} \rightarrow \mathbb{R}^n$ is said to define an {\it admissible nonlinear controller} if the following conditions are satisfied: First,  mapping $K$ is {\it causal}; that is,~$K_0$ is a constant, and, given any~ sample path~$X = (X_0,X_1,...,X_{n-1}) \in {\cal X}$, for~$k = 1,2,...,n-1$, the controller's~$k$-th component~$K_k(X)$ depends only on~$X_0,X_1,...,X_{k-1}$, with the resulting bet size at stage~$k$ given by
$
u_k(K,X) = K_{_{k}}(X)V_k(K,X).
$
Second, controller $K$ should satisfy the {\it budget constraint}~$|u_k(K,X)| \leq V_k(K,X)$, i.e., $|K_{_{k}}(X)|\leq 1$ for all $k=0,\ldots,n-1$. In the sequel, we
denote the set of all admissible controllers by~${\cal K}$.

It is important to note that~${\cal K}$ includes linear controllers as  special case.
More generally, members of ${\cal K}$ can be highly nonlinear functions of~$X$. Accordingly, whenever appropriate,~$K_{k}(X)$ is referred to as a {\it nonlinear feedback gain}.  Now, along sample path~$X \in {\cal X}$, the resulting account value is obtained recursively as $V_0(K,X)=V_0$, and
\[
V_{k+1}(K,X) = (1 +  K_{_{k}}(X)X_k)V_k(K,X),\;\;k=0,1,\ldots
\]
In Figure 1, the binary tree associated with the state transitions above are shown. We also draw attention to the color scheme used for the nodes: At any given stage~$k$, two nodes with the same color represent sample pathes with the same number of heads over the prior stages~$0,1,2,...,k-1$. As seen in the main result to follow, at such nodes, the \emph{optimal} robust nonlinear control, has the same nonlinear gain~$K_k(X)$.
Now, continuing with the analysis, %the final account value is
%$$
%V_n(K,X) = \left(\prod_{k = 0}^{n-1}(1 +  K_{_{k}}(X)X_k)\right)V_0,
%$$
%and, therefore,
the final logarithmic growth is %given by
$$
\frac{1}{n}\log\left(\frac{V_n(K,X)}{V_0}\right)
       = \frac{1}{n} \sum\limits_{k = 0}^{n-1}\log\left((1 +  K_{_{k}}(X)X_k\right).
$$

We now turn our attention to the starting point for much of the analysis to follow: the simple formula for the Expected Logarithmic Growth (ELG) as a function of the controller nonlinear gains~$K \in {\cal K}$ and the probability of heads~$p$. We first find the probability of a sample path $X \in {\cal X}$, given by
$$
P(X) \doteq p^{n_h(X)}(1-p)^{n-n_h(X)}
$$
where~$n_h(X)\doteq\#\{i=0,\ldots,n-1:X_i=1\}$ is the number of heads.
%; i.e.~$ n_h(X)\doteq\#\{i=0,\ldots,n-1:X_i=1\}$.
The ELG as a function of $K$ and $p$ is found as\footnote{Formally, it is possible that $V_n(K,X) = 0$ for some $X\in\mathcal{X}$. If $P(X)>0$ for at least one such sample path, then $\ELG_K(p)=-\infty$;
otherwise, we neglect the resulting summands by using the convention $0\cdot(-\infty) = 0$.
}
\[
\ELG_K(p)=\frac{1}{n}\sum_{X \in {\cal X}}P(X)\sum_{k=0}^{n-1}\log(1 +  K_k(X)X_k).
\]
\begin{figure}[htb]
\centering
\includegraphics[width=1.9in]{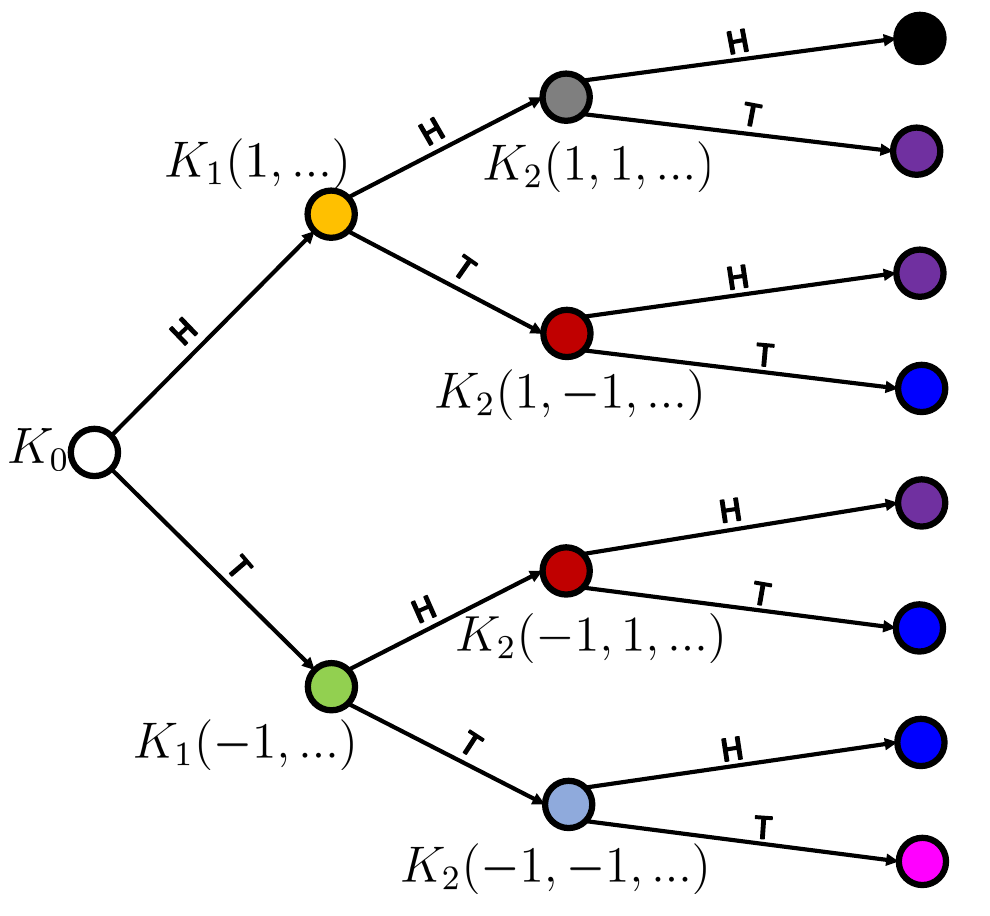}
%\vspace*{-1.5in}
\caption{\bf Causal Controllers and Random Walks on a Binary Tree
%Nodes of same color correspond to equal values of the \emph{optimal} controller $K^*$.
}
\label{fig.tree}
\end{figure}
%\begin{eqnarray*}
%\ELG_K(p)&\doteq&\frac{1}{n}\sum_{X \in {\cal X}} P(X)\log \frac{V_n(K,X)}{V_0}\\
%& = & \frac{1}{n}\sum_{X \in {\cal X}}P(X)\sum_{k=0}^{n-1}\log(1 +  K_k(X)X_k).
%\end{eqnarray*}
%\]

\vskip -.1in
\section{Robustness Formulation}
Per earlier discussion, we now formulate a Robust Expected Logarithmic Growth problem involving uncertainty in the probability of heads. To this end, let~${\cal P} \subseteq [0,1]$ denote a Lebesgue measurable set of the possible values for $p$ against which we seek robustness, and, let~$\mu({\cal P})$ be its corresponding Lebesgue measure. For example, if the only a priori information we have about the probability of heads are bounds
$$
0\leq p_{min} \leq p \leq p_{max} \leq 1,
$$
then with~${\cal P} = [p_{min},p_{max}]$, we have~$\mu({\cal P})= p_{max}-p_{min}$. In the sequel, to avoid trivialities, we assume~$\mu({\cal P})> 0$. Unlike Kelly's perfect-information scenario, we cannot take~$K \in {\cal K}$ depending on the unknown probability of heads~$p$;
the controller $K=K_{\cal P}$ should be determined by the known set~${\cal P}$.

\subsection*{Comparison With Kelly's Perfect-Information Optimum}
To assess the robustness of any particular controller, we compare it's expected logarithmic growth, as a function of $p\in{\cal P}$, with that of Kelly's {\it perfect-information ELG optimum}~$K_p^*\doteq 2p-1$ described in the Introduction.
In this regard, we view~$K_p^*$ as a member the admissible set~${\cal K}$ and a straightforward calculation leads to optimal performance level
$$
\ELG^*(p) = p\log(2p) + (1-p)\log(2(1-p)).
$$
This quantity, serves as our ``gold standard'' against which we assess the robust performance of controllers with {\it imperfect information}. That is, given
any~$K\in {\cal K}$, we first observe that the inequality
$
\ELG^*(p) \geq \ELG_K(p)
$
must hold for all $p\in {\cal P}$ and $K\in\mathcal{K}$. Hence, the associated error integral
$$
Err(K) \doteq \int_{p \in {\cal P }}\left(\ELG^*(p) - \ELG_K(p)\right)dp
$$
is minimized by maximizing the Integral Expected Logarithmic Growth (IELG), that is, the function
$$
\IELG_{K} \doteq  \int_{p \in {\cal P }} \ELG_K(p)dp
$$
over all admissible controllers $K\in\mathcal{K}$.
By viewing~$p$ as a random variable uniformly distributed over $\cal P$, the IELG is in fact proportional to
the \emph{expectation} of random variable $\ELG_K$:
\[
\boldE(\ELG_K) =\frac{1}{\mu(\mathcal{P})}\int_{p \in {\cal P }} \ELG_K(p)dp=\frac{\IELG_K}{\mu(\mathcal{P})}.
\]
It is noted that use of a uniform distribution, implicit in the integral above, is being used solely for the sake of simplicity of the exposition.
As discussed in Section~VII, the analysis to follow is easily modified to address more general distributions on $\mathcal{P}$.
In the theorem to follow in Section~VI, it is seen that element~$K^* \in {\cal K}$ maximizing the IELG exists and is unique, and, we provide a simple and efficient formula to compute it.
\section{The Subclass of Static Linear Controllers}

By way of preliminaries, we say that an admissible~$K \in{\cal K}$ defines an
{\it  admissible static linear controller} if there exists a constant~$K_0\in[-1,1]$ such that
$K(X) \equiv K_0$ for all~$X \in {\cal X}$.
%In the sequel, we use the notation~${\cal K}_0$ to represent this subclass of~${\cal K}$.
In order to study the robust performance over this subclass of controllers, we work with the function
\[
\begin{split}
%f(K_0)&\doteq \int_{p\in\mathcal{P}}\ELG_{{K_0}}(p)dp\\
\IELG_{K_0}=\int_{p\in\mathcal{P}}(p\log(1+K_0)+(1-p)\log(1-K_0))\,dp.
\end{split}
\]
\noindent{\bf Lemma:} {\it The static feedback gain
 maximizing~$\IELG_{K_0}$ subject to the constraint~$K_0\in[-1,1]$ is unique and given by}
$$
K_0^* =2\bar p-1;\quad \bar p\doteq\frac{1}{\mu({\cal P})}\int_{p \in {\cal P }}pdp.
$$
\noindent {\bf Proof}: A straightforward computation indicates that for $K_0\in(-1,1)$, function~$f(K_0)\doteq\IELG_{K_0}$ has the first derivative
\[
f'(K_0) = \frac{2\int_{\cal P}pdp-\mu(\mathcal{P})-K_0\mu(\mathcal{P})}{1-K_0^2}=\frac{\mu(P)(2\bar p-1-K_0)}{1-K_0^2},
\]
whereas $f(\pm 1)=-\infty$. Since this derivative is positive when $K_0<K_0^*$, negative when $K_0>K_0^*$ and zero at~$K_0=K_0^*$, it follows that~$K_0^*$
maximizes $\IELG_{K_0}$ over $K_0\in[-1,1]$. $\square$

\subsubsection*{Observations and Important Special Cases}
Notice that the optimal robust static linear control coincides with the Kelly's formula~$K^*_p$ for the perfect-information case with probability
of heads~$p$ being the \emph{centroid} $\bar p$ of uncertainty set $\mathcal{P}$. It is also interesting that we can obtain Kelly's formula~$K^*_p$ as a special case of our robustness analysis by relaxing the positive measure assumption on~${\cal P}$ and using a limiting argument. That is, suppose~$p \in (0,1)$ and~${\cal P} = [p, p +\delta]$ with parameter~$\delta$ satisfying~$0 < \delta \leq 1-p$. Then, taking the limit as~$\delta \rightarrow 0$, it is straightforward to verify that the optimal robust linear feedback gain reduces to~$K_0^*= 2p-1$. Next, we consider the case when~${\cal P}$ is a union positive-length disjoint intervals
$$
{\cal P} = \bigcup\nolimits_{{i =1}}^m[p_{min,i},p_{max,i}].
$$
Applying the lemma, a straightforward calculation leads to
$$
K_0^* = \frac{\sum_{i=1}^m ({p_{max,i}^2 - p_{min,i}^2})}
          {\sum_{i=1}^m ({p_{max,i} - p_{min,i}})}-1
$$
which specializes further: If all differences~$p_{max,i} - p_{min,i}$ are the same, then, the formula above reduces to
$$
K_0^* = \frac{1}{m} \sum\nolimits_{i =1}^m (p_{min,i} + p_{max,i}) - 1.
$$
which, for the single-interval case~${\cal P} = [p_{min},p_{max}]$ becomes
$$
K_0^* = p_{min} + p_{max} -1.
$$
\section{Examples: Optimal Nonlinear Controller}
To motivate the key ideas underlying the general result in the theorem to follow, we calculate the optimal robust nonlinear control for the simpler special cases~$n =2$ and~$n = 3$. As seen below, each of these optima can be found in a very simple manner.
That is, each of the desired nonlinear feedback gains comprising the optimum~$K^* \in {\cal K}$ is found via single-variable maximization whose solution
admits a closed form. The simple argument used in these examples is also important in the proof for the general case of~$n >2$ in the theorem to follow.

\subsubsection*{Example~1} Beginning with the $n=2$ and ${\cal P} = [0,1]$, to simplify calculations, we first represent the nonlinear controller components~$K_k(X)$ employing the shorthand notation~$a = K_0(X)$ for all~$X \in {\cal X}$, $b = K_1(1,1) = K_1(1,-1)$ and~$c = K_1(-1,1) = K(-1,-1)$, we first calculate
\[
\begin{aligned}
\ELG_K&(p)=\\
        =\tfrac{1}{2}&\left[p^2\log(1+a)(1+b)+ p(1-p)\log(1+a)(1-b)+\right.\\
         p(1&-p)\log(1-a)(1+c)+(1-p)^2\log(1-a)(1-c)\left.\right].
\end{aligned}
\]
Then, upon expanding the logarithms above and integrating with respect to~$p\in{\cal P}$, we obtain the IELG function
\[
\begin{aligned}
\IELG_K  = & \;\;\;\tfrac{1}{2}\left[ \tfrac{1}{2}\log(1+a) + \tfrac{1}{2}\log(1-a)+\right.\\
         & \;\;\;\; +\;\tfrac{1}{3}\log(1+b) + \tfrac{1}{6}\log(1-b)+\\
         & \;\;\;\left. +\;\tfrac{1}{6}\log(1+c) + \tfrac{1}{3}\log(1-c)\right].
\end{aligned}
\]
Next, we note that the desired optimization with respect $(a,b,c)$ above can be solved by maximizing three separate single-variable strictly concave functions; each of
functions is of the form~$f(x) = \alpha\log(1+x) + \beta\log(1-x)$ with unconstrained maximum, obtained by setting the derivative to zero, given by~$x^* = (\alpha - \beta)/(\alpha + \beta)$. Then, we simply observe that the budget-constrained optimum is also~$x^*$ because in all three cases above,~$\alpha$ and~$\beta$ are such that~$|x^*| \leq 1$. Based on these considerations, we immediately arrive at a unique maximizing solution~$a = 0; b = \frac{1}{3};c =-\frac{1}{3}$. In other words, it is optimal to skip the first bet, and  bet $1/3$ of the account value on heads if the first toss comes up to be heads ($X_0=1$), and bet $1/3$ of the account value on tails if the first toss comes up to be tails ($X_0=-1$). Furthermore,
via a straightforward substitution we obtain~$\ELG_{K^*}= \frac{1}{3}\log\frac{32}{27}\approx 0.0566.$ Note that the optimal static \emph{linear} controller from Lemma~1 is degenerate: $K_0^*=0$ and $\ELG_{K_0^*}=0$ (no bets are made).

\subsubsection*{Example~2} Now proceeding to the analysis for~$n = 3$, we consider the case~${\cal P} = [0.25,0.95]$ and again represent the nonlinear control gains using a shorthand notation which takes causality into account: For each~$K_k(X)$, only the first~$k$ components of its argument~$X$ are indicated and we take
\[
\begin{gathered}
a=K_0;\;b=K_1(1);\;c=K_2(1,1);\;d=K_2(1,-1);\\
e=K_1(-1);\;f =K_2(-1,1);\;g=K_2(-1,-1).
\end{gathered}
\]

Now, we repeat the straightforward sum-of-logarithms computation as in Example~1 and integrate~$\ELG_K(p)$ to arrive at
\[\begin{aligned}
\IELG_K  = &\;\tfrac{1}{3}\left[0.09333\log(1 - a) + 0.1400\log(1+a)\right.\\
     & +0.04647\log(1-b) + 0.09353\log(1+b)\\
   &+ 0.02598\log(1-c) + 0.06755\log(1+c)\\
   & + 0.02049\log(1-d)+ 0.02598\log(1+d)\\
  & + 0.04686\log(1-e) + 0.04647\log(1+e)\\
 & + 0.02049\log(1-f)+ 0.02598\log(1+f)\\
  & \left.+ 0.02637\log(1 - g) + 0.02049\log(1+g)\right].
\end{aligned}\]
Now, proceeding again as in Example~1, we maximize separately with respect
to each of the seven control parameters and obtain unique optimum~$a = 0.2$, $b = 0.3361$, $c = 0.4445$, $d = 0.118$, $e = -0.004167$, $f = 0.118$, $g = -0.007429$. Finally, for performance comparison purposes, we calculate the optimal static gain~$K_0^* = p_{min} + p_{max} -1 = 0.2$.

It is interesting to compare the robust performance for the optimal nonlinear controller~$\ELG_{K^*}(p)$ with that of the optimal static linear controller~$\ELG_{K_0^*}(p)$, as functions of the probability~$p \in {\cal P}$ and benchmarked with Kelly's perfect-information optimum~$\ELG^*(p)$ as the best-possible upper bound. In Figure~2 below, where these quantities are plotted, the following
is noted: While the optimal nonlinear controller leads to a larger IELG (the area under the curve), at some specific values of $p$ the static linear controller outperforms it in terms of $\ELG(p)$. This is not surprising, because at the midpoint $p=0.6$ of interval for $\mathcal{P}$, the optimal static gain $K_0^*$, in fact, coincides with the Kelly gain $K^*(p)$, which maximizes the $ELG_K(p)$ among \emph{all} admissible controllers $K\in\mathcal{K}$.
\begin{figure}[tb]
\centering
\includegraphics[width=0.75\columnwidth]{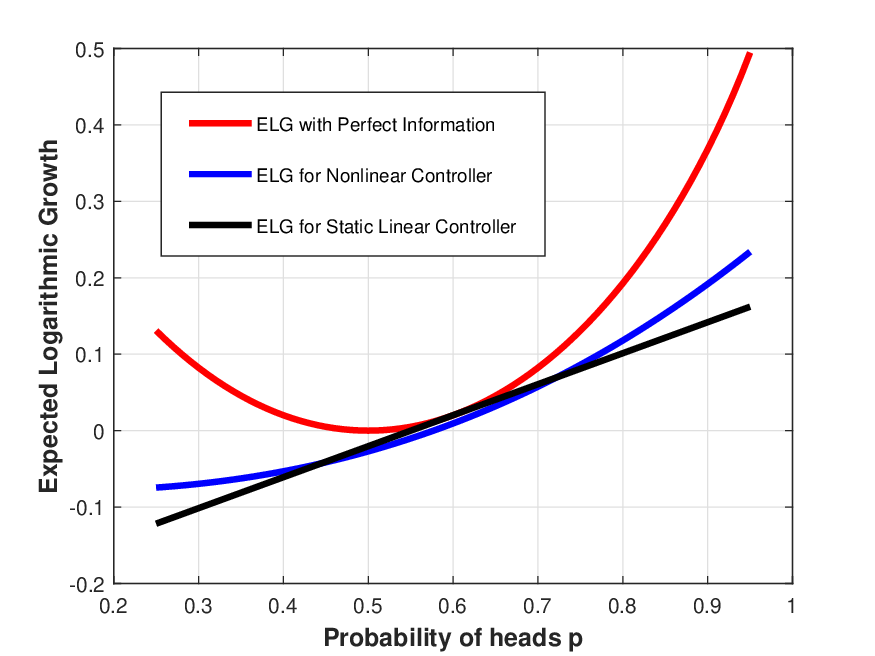}
%\vspace*{-1.5in}
\caption{\bf Robust Performance Plots for Comparison Purposes}
\end{figure}

\section{Main Result on Robust Optimal Control}

The theorem below, establishes the existence and uniqueness of an optimal robust nonlinear controller~$K^* \in {\cal K}$ and characterizes it with an explicit formula. Although there are~$2^n-1$ nonlinear controllers gains associated with the nodes, the theorem tells us that at each stage~$k$, there are only~$k+1$ possible values of for the optimal robust nonlinear gain~$K_k^*(X)$. Summing up these numbers across all stages, we see that  the total number of nonlinear gains to be calculated is~$1+2+\ldots+n=n(n+1)/2$; i.e., the computational burden of finding all gains $K_k(X)$ increases quadratically in~$n$ rather than exponentially. On the other hand, if the time between consecutive coin flips is suitably large, the controller in the theorem below can be implemented ``dynamically'' with no need to pre-compute the optimal nonlinear gains~$K_k^*(\bar X)$.\\
\\
%\vskip1mm
\noindent
{\bf Theorem}: {\it The integral expected logarithmic growth $\IELG_K$
%$$
%f(K) \doteq  \int_{p \in {\cal P }} \ELG_K(p)dp,
%$$
defined over on the set of admissible controllers $K\in\mathcal{K}$, has a unique maximizer~$K^*$, whose nonlinear control gain at stage~$k$ for a  sample path $\bar X\in\mathcal{X}$, is given by
\[
\begin{gathered}
K_k^*(\bar X)=\frac{\int_{p\in\mathcal{P}}p^{q_k}(1-p)^{k-q_k}(2p-1)\,dp}
{\int_{p\in\mathcal{P}}p^{q_k}(1-p)^{k-q_k}\,dp},%,\\
%q=q_k(X)\doteq\#\{i=0,1,\ldots,k-1:X_i=1\}
\end{gathered}
\]
where $q_k=q_k(\bar X)=\#\{i=0,\ldots,k-1:\bar X_i=1\}\leq k$ is the number of heads occurring over the first $k$ coin flips}\footnote{By definition, $q_0(\bar X)=0$, i.e., the optimal control gain $K^*_0(X)\equiv K^*_0$ at stage $k=0$ coincides with the optimal linear controller from Lemma~1.}.
\\
%\bb{some of the steps in the proof below require attention to avoid ill-defined
%sets, empty sets, quantities etc for $k= 0$ and~$n = 1$}

\subsubsection*{Proof}\label{sec.proof}
In the arguments to follow, sample path $\bar X \in {\cal X}$ and stage number~$k \in \{0,1,...,n-1\}$ are assumed to be fixed, and we
let ${\cal X}_k(\bar X) \doteq {\cal X}_k^+(\bar X)\cup {\cal X}^-_k(\bar X)$, defining sets ${\cal X}_k^{\pm}$ by\footnote{For $k=0$, sets~${\cal X}_k^+$ and ${\cal X}_k^-$ consist of all sample paths starting from $X_0=1$ and $X_0=-1$, respectively.}
\[
\begin{gathered}
{\cal X}_k^+(\bar X) \doteq \{X \in {\cal X}:
 X_i = \bar X_i\; \forall i = 0,1,...,k-1;\, X_k = 1\},\\
{\cal X}_k^-(\bar X) \doteq \{X \in {\cal X}:
 X_i = \bar X_i\; \forall i = 0,1,...,k-1;\, X_k = -1\},
\end{gathered}
\]
%wih sets~${\cal X}_0^+$ and ${\cal X}_0^-$ consisting of all sample paths starting from $X_0=1$ and $X_0=-1$, respectively.
the expected logarithmic growth function can be written as
\begin{equation*}
\begin{aligned}
\ELG_K(p)& = &\frac{1}{n}\sum_{X \in {\cal X}_k^+(\bar X)}\hskip -.15in P(X)\sum_{i = 0}^{n-1}
\log(1 + K_i(X)X_i)\\
& & \;\;+\frac{1}{n}\sum_{X \in {\cal X}_k^-(\bar X)}\hskip -.15in P(X)\sum_{i = 0}^{n-1}\log(1 + K_i(X)X_i)\\
& & \;\;+\frac{1}{n}\sum_{X \notin {\cal X}_k(\bar X)}
\hskip -.15in P(X)\sum_{i = 0}^{n-1}\log(1 + K_i(X)X_i).
\end{aligned}
\end{equation*}
Next we note that all of the terms above involving~$K_k(\bar X)$ can be isolated by setting~$i = k$ in the first two terms above. As far as the third term is concerned, it is independent of~$K_k(\bar X)$ because  admissibility of the controller forces~$K_k(X)= K_k(\bar X)$ for all~$X\in {\cal X}_k(\bar X)$. Now integrating $\ELG_K(p)$ over $p\in\cal P$, it is straightforward to see that maximization of~$\IELG_K$ over $K\in{\cal K}$ reduces to maximization of
the single-variable function
$$
g_{k,\bar X}(K_k)\doteq\alpha_k(\bar X)\log(1 + K_k)+ \beta_k(\bar X) \log(1 - K_k)
$$
over the interval $K_k\in[-1,1]$. Here $\alpha_k,\beta_k$ are defined as
\[
\begin{gathered}
%g_k(\bar X,K_k)\doteq\alpha_k(\bar X)\log(1 + K_k)+ \beta_k(\bar X) \log(1 - K_k),\\
\alpha_k(\bar X)\doteq  \int_{p\in\cal P }
 \sum_{X \in {\cal X}_k^+(\bar X)}\hskip-.15in P(X)\,dp=\int_{p\in\cal P}\boldP(X\in\mathcal{X}_k^+(\bar X))\,dp,\\
 \beta_k(\bar X)\doteq\int_{p\in{\cal P }}
 \sum_{X \in {\cal X}_k^-(\bar X)} \hskip -.15in P(X)\,dp=\int_{p\in\cal P}\boldP(X\in\mathcal{X}_k^-(\bar X))\,dp.
\end{gathered}
\]
%function $f$ is found as
%\[
%f(K)=\frac{1}{n}\sum_{k=0}^{n-1}\sum_{\bar X\in\cal X}g_k(\bar X,K_k(\bar X));
%\]
%where the function $g_k(\bar X,K_k)$ is readily obtained as
%
%are defined as follows
%\[
%\begin{gathered}
%g_k(\bar X,K_k)\doteq\alpha_k(\bar X)\log(1 + K_k)+ \beta_k(\bar X) \log(1 - K_k),\\
%\alpha_k(\bar X)\doteq  \int_{p\in\cal P }
% \sum_{X \in {\cal X}_k^+(\bar X)}\hskip-.15in P(X)\,dp=\int_{p\in\cal P}P(X\in\mathcal{X}_k^+(\bar X))\,dp,\\
% \beta_k(\bar X)\doteq\int_{p\in{\cal P }}
% \sum_{X \in {\cal X}_k^-(\bar X)} \hskip -.15in P(X)\,dp=\int_{p\in\cal P}P(X\in\mathcal{X}_k^-(\bar X))\,dp.
%\end{gathered}
%\]
To find $\alpha_k$, notice that event $X\in\mathcal{X}_k^{+}$ is the intersection of
$k+1$ events $X_i=\bar X_i$ (where $i=0,\ldots,k-1$) and $X_k=1$; these events are mutually independent.
Among these events, there are $1+q_k(\bar X)$ events of type $X_j=1$ and probability $p$ and $k-q_k(\bar X)$ events of type
$X_j=-1$ and probability $1-p$. Multiplying these probabilities, one has
\[
\alpha_k(\bar X)=\int_{p \in {\cal P }}p^{q_k(\bar X)+1}(1-p)^{k-q_k(\bar X)}dp.
\]
A very similar argument yields in the expression for $\beta_k$; i.e.,
%
% $X\in\mathcal{X}_k^{-}$ is the intersection of $k+1$ independent events $X_i=\bar X_i$ (where $i=0,\ldots,k-1$) and $X_k=-1$, among which
%$q_k(\bar X)$ events have probability $p$ and $k+1-q_k(\bar X)$ events have probability $1-p$. Hence,
\[
\beta_k(\bar X)=\int_{p \in {\cal P }}p^{q_k(\bar X)}(1-p)^{k+1-q_k(\bar X)}dp.
\]
%Now, to complete the proof, we note that
%A controller $K^*$ maximizes $f(K)$ over all $K\in\mathcal{K}$ if and only if, for each sample path $\bar X$, $K_k^*(\bar X)$ is a maximizer
%of function $g(\bar X,K_k)$ over all $K_k\in[-1,1]$.
Since $g_{k,\bar X}(K_k)$ is strictly concave on interval $K_k\in[-1,1]$ and $g_{k,\bar X}(\pm 1)=-\infty$, the optimum is found by setting the derivative with respect to~$K_k$ to zero; i.e., we obtain
%$%\partial g_k(\bar X,K_k)/\partial K_k$ to zero, that is,
\[
\begin{aligned}
K_k^*(\bar X) &= \frac{\alpha_k(\bar X)-\beta_k(\bar X)}{\alpha_k(\bar X)+\beta_k(\bar X)}=\\
&=\frac{\int_{p\in\mathcal{P}}p^{q_k(\bar X)}(1-p)^{k-q_k(\bar X)}(2p-1)\,dp}
{\int_{p\in\mathcal{P}}p^{q_k(\bar X)}(1-p)^{k-q_k(\bar X)}\,dp},
\end{aligned}
\]
which satisfies $|K_k^*(\bar X)|\leq 1$ as required. $\square$
\subsection*{Nonlinear Versus Linear Control}
We are now prepared to address one of our main contentions articulated the title and abstract. That is, except for the trivial case of single-flip game~($n = 1$), we establish, as a corollary of the theorem, that the optimal nonlinear controller~$K^* \in {\cal K}$ robustly outperforms the optimal static linear feedback~$K_0^*$.
\vskip1mm
\noindent
{\bf Corollary}: { \it For $n > 1$ steps, one has $\IELG_{K^*}>\IELG_{K_0^*}$.}
\vskip2mm
To facilitate the proof,  we first provide a preliminary lemma.
\vskip1mm\noindent
{\bf Preliminary Lemma}: {\it Let ${\cal P} \subseteq [0,1]$ be a Lebesgue measurable set with~$\mu({\cal P}) > 0$. Then, for all $n>1$, it follows that }
$$
\int_{p \in{\cal P}} p^{n-1}\,dp \int_{p \in {\cal P}} p\,dp < \mu({\cal P})\int_{p \in {\cal P}} p^n\,dp.
$$
\noindent {\bf Proof.} Denoting $I_k\doteq\int_{p \in{\cal P}} p^k\,dp$, we need to prove that $I_{n-1}I_1<I_0I_n$. %We proceed in three steps.
%\\
Applying H\"older's inequality to pairs of functions $f(p)=p,g(p)=1$ and $g(p)=1,h(p)=p^{n-1}$
and two conjugate exponents $n$ and $m=n/(n-1)$, one obtains
\[
\begin{gathered}
I_1=\|fg\|_{L_1(\cal P)}<\|f\|_{L_n(\cal P)}\|g\|_{L_{m}(\cal P)}=I_n^{1/n}I_0^{1/m},\\
I_{n-1}=\|gh\|_{L_1(\cal P)}<\|g\|_{L_{n}(\cal P)}\|h\|_{L_m(\cal P)}=I_0^{1/n}I_n^{1/m}.
\end{gathered}
\]
Now multiplying these two inequalities, it follows that $I_1I_{n-1}<I_n^{1/n}I_0^{1/m}I_{n}^{1/m}I_0^{1/n}=I_nI_0$ finishing the proof. $\square$

\subsubsection*{Proof of Corollary} Recalling the theorem, $K^*$ is the \emph{unique maximizer} of  $IELG_K$ over all admissible nonlinear gains~${K \in \cal K}$. Since the optimal linear static gain $K_0^*$ is admissible, it suffices to show that $K^*_k(X)\ne K^*_0$ for at least one $k \in \{0,1,\ldots,n-1\}$ and at least one sample path $X \in {\cal X}$. Indeed, considering the distinguished sample path corresponding to all heads; i.e., $X_0=\ldots=X_{n-1}=1$, we first note that~$q_{n-1}(X)=n$. Recall that $K^*_0=2\bar p-1$, where $\bar p$ is the centroid of $\mathcal{P}$. Now applying the theorem, we have
\[
K^*_{n-1}(X)=2\hat p_n-1;\;\;\hat p_n\doteq \frac{\int_{p\in\mathcal{P}}p^{n}\,dp}{\int_{p\in\mathcal{P}}p^{n-1}\,dp}>\bar p
\]
with the latter inequality implied by our Preliminary Lemma. Thus, $K^*_{n-1}(X)>K_0^*$, and hence $K_0^*$ is non-optimal. $\square$

\section{Discussion of Two Generalizations}
As previously mentioned, the assumption of equal payoffs for heads and tails was made solely for simplicity of exposition and brevity of the presentation. We now sketch the key ideas indicating how the robustly optimal nonlinear gains~$K_k^*(X)$ are obtained for the unequal payoff case with~$X_k= a$ for heads and $X_k= b$ for tails at stage $k$, and, to avoid trivialities, it is assumed that~$ b < 0 < a$.   Indeed, let $q_{a,k} \doteq q_{a,k}(X)$ be the number of heads seen on the first~$k$ flips and
\[
\begin{gathered}
\alpha_{_{a,k}} \doteq \int_{{\cal P}}p^{q_{_{a,k}}+1}\left(1-p\right)^{k - q_{_{a,k}}} dp,\\
\beta_{_{a,k}} \doteq \int_{{\cal P}}p^{q_{_{a,k}}}\left(1-p\right)^{k - q_{_{a,k}}+1} dp.
\end{gathered}
\]
Then, to get the robustly optimal gain~$K_k^*=K_k^*(X)$, we form the \emph{strictly concave} scalar function
$$
g_k(K_k) \doteq \alpha_{_{a,k}}\log(1+ K_ka)+ \beta_{_{a,k}}\log(1+ K_kb)
$$
to be maximized subject to budget constraint $-1 \leq K \leq 1$ and the
requirement~$V_{k+1} \geq 0$ associated with both well-definedness of the logarithms above and bankruptcy considerations. Accordingly, with
$$
m \doteq \max\left\{-1,-\frac{1}{a}\right\};\;\;  M  \doteq \min\left\{1,\frac{1}{|b|}\right\}
$$
we obtain the optimal nonlinear gain~$K_k^*$ as the unique maximizer of~$g_k(K_k)$ on the interval~$
K_k \in [m,M]$ and observe that our main Theorem corresponds to the case where $a=M=1$ and $b=m=-1$.

Finally, as has been already mentioned earlier in the paper, our results retain their validity if one replaces the uniform distribution on $\cal P$ by a positive finite measure $\hat\mu$ (defined, at least, on Borel subsets of $\cal P$) and redefining the IELG as
\[
\IELG_K=\int_{\cal P}\ELG_K(p)\hat\mu(dp).
\]
Our Lemma, Theorem and Corollary then remain valid, replacing $dp$ in all integrals by $\hat\mu(dp)$ and $\mu(P)$ by $\hat\mu(P)$.

\section{Conclusion and Future Research}

In this paper, our main objective was to demonstrate that nonlinear control has an important role to play in large classes of betting games dealing with Expected Logarithmic Growth.  To this end, we considered a simple coin-flipping game as a demonstration case to convey our main ideas. Whereas a static linear control is ``unbeatable'' with a perfectly known probability of heads~$p$, this does not hold true when robustness with respect to variations in~$p$ is of concern. For this situation, we showed that the optimal controller with its nonlinear gains~$K^*$ robustly outperforms the optimal static linear controller with its gain~$K_0^*$.

Perhaps the main implication of our results is that future
study of nonlinear control with more general problem formulations, relevant for the field of finance, is likely to bear fruit.
By way of future research, in addition to the generalizations sketched in Section~VII, we believe that it should be possible to address the case when the returns~$X_k$ take on multiple or even a continuum of values governed by rather general probability distributions. Another possible generalization begins with ``vector sample paths'' $X$ in lieu of the scalar ones considered here. Such a formulation can be viewed in a robust portfolio balancing context with results along the lines serving  as a stepping stone to applications such as algorithmic stock trading in financial markets. In bringing such results from theory to practice, it would be important to add terms to the account value dynamics for~$V_k$ to
include consideration of factors such as the risk-free and margin interest rates, leverage and transaction costs.

Finally, we mention one additional continuation of this research which is motivated by the following observation: An adaptive controller aimed at maximizing expected logarithmic growth, say along the lines of those given in recent papers such as~\cite{Despons_et_al_2022} and \cite{Dettu_et_al_2022}, should rightfully be viewed as a member of our admissible control set~${\cal K}$. Accordingly, our plan for future research involves exploring the connection between results in adaptive and nonlinear control which have traditionally been viewed as rather separate areas. In this regard, further motivation for such work is provided by the simple example provided for ~$n = 2$. For this low-dimensional example, our optimal three-gain robust nonlinear controller turns out to be the same as the one provided in~\cite{Despons_et_al_2022}.
It should be also noted that the ``adaptive Kelly'' gain~\cite{Despons_et_al_2022}
$2\hat p_k(X)-1$ converges, as $k\to\infty$, to the ideal Kelly gain $K^*(p)=2p-1$ with probability $1$ due the Law of Large Numbers.
Here $\hat p_k(X)$ is the estimated probability of heads inferred from sample path $X$ (e.g., constructed as
in~\cite{Despons_et_al_2022}).
Since the robust optimal controller provides the value of IELG that is not less than the IELG of ``adaptive Kelly'' control gain
(being admissible), it can be proved that the optimal IELG converges, as $n\to\infty$, to the IELG of ideal Kelly controller.
We leave the rigorous analysis (with convergence rate estimates) for future research.
\vspace{2mm}
\normalfont

\end{document}